\documentclass[a4paper]{article}
\usepackage{amsfonts}
\usepackage{times}
\usepackage{mathptm}
\title{Lagrangian Matroids Associated\\ with Maps on Orientable Surfaces
}
\author{Richard F. Booth \\ Department of Mathematics \\ UMIST \\ PO Box 
88 \\ Manchester M60 1QD \\ United Kingdom
\and 
Alexandre V. Borovik \\ Department of Mathematics \\ UMIST \\ 
PO Box 88 \\ Manchester M60 1QD \\ United Kingdom
\and Israel Gelfand \\ Department of Mathematics \\ Rutgers 
University \\ New Brunswick NJ 08903 \\ USA}
\date{28 January 1999\footnote{updated 2000/10/18 for the web}}
\newcommand{\qed}{\ \hspace*{\fill} $\diamond$ \vskip 5pt}
\newenvironment{proof}{\paragraph{Proof}}{\qed}
\newtheorem{lem}{Lemma}
\newtheorem{thm}[lem]{Theorem}

\newcommand{\B}{{{\cal B}}}

\newcommand{\M}{{{\cal M}}}
\newcommand{\Fl}{{\overline{F}}}

\begin{document}

\pagestyle{myheadings}
\markright{{\small R.~F.~Booth et al., 
Lagrangian Matroids and Maps on  Orientable Surfaces}}
\maketitle

\section*{Introduction}

The aim of this paper is to clarify the nature of combinatorial 
structures associated with 
maps on closed compact surfaces. These are fairly classical objects; 
however, not long  
ago it was discovered by A.~Bouchet that maps are associated with
$\Delta$-matroids \cite{bou:map} (or Lagrangian matroids in the  terminology 
of \cite{bgw}).  $\Delta$-matroids are related to maps in almost 
the same way as ordinary matroids to graphs. 

In this paper we explore this parallel further and show that the 
resulting Lagrangian matroids are representable by matrices in 
a sense analogical to representability of ordinary matroids, 
thus tranferring the classical Rado's theorem \cite{rado} to 
Lagrangian matroids. The proper setting for the representation 
is provided by cohomology of the surface.
Our proof is very elementary. It is worth mentioning that
David Stone \cite{stone} found a cohomological proof
of the main technical result of the paper, Theorem~\ref{cyciso}.
His proof is not
much shorter than ours but promises a possible
generalisation of our result to higher dimension.

Our next observation is that the classical `spanning tree algorithm' 
for graphs (which becomes the `greedy algorithm' of matroid theory 
and the theory of $\Delta$-matroids \cite{bou}) can be interpreted in 
our setting  in a most elementary way, as a `peeling off' procedure 
which cuts the (connected) surface into one closed ring-shaped 
peel, or, in some degenerate cases, to a $2$-cell peel.
This procedure is local, that is,
at every step uses information only about a small part of the surface 
around the knife. 



\section{Orthogonal matroids and maps}

\subsection{Maps on compact surfaces}

A {\em graph} is a compact space $G$ partitioned into a finite set 
$E=E(G)$ of open $1$-cells called {\em edges} and a
finite set $V=V(G)$ of points called {\em vertices} such that every 
edge has a vertex at each end (possibly the same vertex at
both ends, in which case we say that the edge is a loop). 
If $S$ is a compact surface and
$G \subseteq S$ we say that the graph is  drawn on the 
surface $S$. 
We shall be concerned only with orientable
connected compact surfaces in this paper.
If connected components of $S\backslash G$ are open $2$-cells,
then we say that $\M := (G,S)$ is a {\em map}, and
the the connected components of $S\setminus G$ are
called the {\em faces} or {\em countries} of the map. 
Thus a map $\M$ introduces on $S$ a structure of a  
$2$-dimensional cell complex.  We assume that the set of edges
(correspondingly, coedges) of $\M$ is non-empty, 
excluding from consideration a trival  map on a sphere
with one vertex and no edges. 


A dual map $\M^*=(G^*,S)$ is defined in the usual way. Its vertices  
are in one-to-one correspondence with the faces of $\M$; 
for each face of $\M$, we choose a point inside it as
the corresponding vertex of $\M^*$. The edges are
chosen such that every edge of $e$ of $G$ is disjoint from every edge 
of $G^*$ except for 
its dual edge, $e^*$, which has with $e$ exactly one point in common. 
We call the $e^*$ 
the coedges. Then it can be shown that the faces of 
$\M^*$ correspond to the vertices of $\M$: every face of 
$\M^*$ contains a unique vertex of $\M$, and every vertex of $\M$
belongs to a unique face of $\M^*$.  
The vertices of $\M^*$ are called {\em covertices}
and the set of covertices is denoted by $V^*$.
Finally, notice that $\M$ is a dual map to $\M^*$.

For $F \subseteq E$ we define $F^*:= \{\,e^* \mid e\in F \,\}$ and
$\Fl := E\backslash F$. We call a set $F\subset E \cup E^*$ 
{\em admissible} if $F\cap F^* = \emptyset$.
It will be convenient to index edges in $E$ by elements of
$I= \{\,1,\ldots,n\,\}$ and the coedges in $E^*$ by the corresponding
elements of $I^*$, so that
$$
E = \{\, e_1, \ldots, e_n\,\} \hbox{ and }
E^* = \{e_{1^*},\ldots, e_{n^*}\,\}.
$$

\subsection{Symplectic and orthogonal matroids}

For a fuller exposition of the theory and definitions of symplectic 
matroids, see  \cite{bgw} and the forthcoming book \cite{book}. Let
$$I=\{1,2,\ldots,n\} \hbox{ and } I^*=\{1^*,2^*,\ldots,n^*\}$$
and $J=I \cup I^*$. Define maps $^*:I \rightarrow I^*$ by $i 
\mapsto i^*$
and
$^*:I^* \rightarrow I$ by $i^*\mapsto i$, so that $^*$ is an 
involutive
permutation of $J$.
Then we say that a subset $K \subset J$ is {\em admissible} if and 
only if
$K \cap K^* = \emptyset$. We denote by $J_k$ the set of all 
admissible $k$-subsets of $J$.

Let $\B \subseteq J_n$ be a a set of admissible $n$-subsets of $J$ 
and let
$M$ be the triple
$(J,\, ^*, \B)$. Then $M$ is a {\em Lagrangian matroid} if it 
satisfies the  {\bf Symmetric Exchange Axiom}:
\begin{quote}
 {\em  For any $A,B \in \B$ and $k\in A\triangle B$ there exists some $i\in 
A\triangle B$ such that
 $A\triangle \{k,k^*,i,i^*\} \in \B$.}
\end{quote}

Here $\triangle$ is the symmetric difference;
$A\triangle B := (A\cup B) \backslash (A\cap B)$.
This axiom is due to Bouchet \cite{bou}, Dress and Havel 
\cite{bdh,dress-havel}, where Lagrangian symplectic matroids appeared,
cryptomorphically, under the  names of
{\em $\Delta$-matroids\/} or {\em metroids}.
We call $\B$ the set of {\em bases} of $M$.
A Lagrangian matroid is a special case of a symplectic matroid; in a 
general symplectic matroid, the bases are elements
of $J_k$ for some $k$. An appropriate axiom system is given in 
\cite{bgw}. In this paper, we shall only be concerned with
Lagrangian matroids.
An {\em orthogonal matroid\/} or {\em even matroid\/}
is a symplectic matroid  
in which the difference between
the number of starred elements in any two bases  is
always even.

We say that an admissible set $F\subset E\cup E^*$
is  {\em independent\/}  if $S\backslash cl(F)$ is connected.
A {\em basis\/} of $M$ is a maximal independent admissible set.

Bouchet \cite{bou:map} proved 
that the set
${\cal B} = {\cal B}(M)$ of all bases of a  map
$M$ is a  even Lagrangian matroid on the set $E \cup E^*$.
This paper provides another proof of Bouchet's result with a
important improvement of {\em representability};
necessary terms are explained in Subsection~\ref{subsec:representability}.

\begin{thm} Let\/ $\M$ be a map with\/ $n$ edges 
on a orientable compact closed surface $S$. Then 
\begin{itemize}
\item 
all bases of\/ $\M$ have  cardinality $n$ and 
\item
the set\/ $\cal B$ of all bases is an orthogonally representable
over\/ $\mathbb{Q}$,  orthogonal Lagrangian matroid.
\end{itemize}
\label{th:map=Lagrange}
\end{thm}

It is well-known that (ordinary) matroids associated
with graphs are representable by matrices 
(\cite{rado}, the exposition can be found in
any book on matroid theory, see, for example, \cite{welsh}).
Our Theorem~\ref{th:map=Lagrange} is a generalisation of this result. 
We shall see in the next subsection that the classical 
Spanning Tree Algorithm for graphs, which corresponds to the Greedy
Algorithm of matroid theory, also has a natural analogue for maps.

\subsection{The greedy algorithm and peeling the skin.}


Draw on each face of $M$ the segments connecting the
covertex with the vertices of this face,
we shall call them {\em diagonals}.
Edges, coedges and diagonals define a baricentric subdivision of $M$.
Each of the triangles of this subdivision has one side a diagonal, one side
 a half-edge of the map, and one side a half-coedge.
It can be immediately seen that, when cutting along
the edges and coedges of a basis,
we obtain a ring (topologically, a punctured 2-cell) or a $2$-cell.
This is due to the fact
that, in view of Theorem~\ref{th:map=Lagrange}, 
the basis forms an admissible set of $n$ elements, and thus every
triangle of the baricentric subdivision is cut along precisely one edge.
This means that each triangle has two neighbors in the modified surface;
hence it is a ring or a $2$-cell. 

The procedure for peeling is the following. 
We construct a sequence of triangles $S_0, S_1, \ldots, S_n$
in the baricentric
subdivision in which $S_i$ and $S_{i-1}$, $i = 1,2,\ldots$,
are adjacent, i.e.\  have in common an edge, coedge or diagonal.

\bigskip
\noindent
{\sc Peeling off Procedure}

\begin{description}
\item[$0^\circ$]  Start at arbitrary triangle $S_0$ in the subdivision. 

\item[$1^\circ$]  Assume that we have constructed $S_0,S_1, \ldots, S_i$. 
\begin{itemize}
\item  If  neither edge $e$  nor coedge  $e^*$ which bound
$S_i$ has  been cut along at the previous steps of the procedure,
cut the surface along  the  entire edge $e$  or coedge $e^*$
chosing $e$ or $e^*$ so that the surface remains connected.
Take for $S_{i+1}$ the triangle which has a common with $S_i$   
coedge (correspondingly, edge) which was  not cut.

\item If  either edge $e$  or coedge  $e^*$  has been  cut
along at one of the previous steps of the procedure, 
take for $S_{i+1}$ the  triangle 
 which lies across the diagonal or noncut
 (co)edge from $S_i$ and is different from
$S_{i-1}$.
\end{itemize}

\item[$2^\circ$] If step $1^\circ$ cannot be made or
if $S_{i+1}$ defined by rule $1^\circ$ is  
one of the triangles $S_0,S_1,\ldots, S_{i-1}$, stop.

\end{description}

We may reformulate this description in terms of 
independent sets. By definition, 
an independent set  of  (co)edges 
is characterised by the property that it is admissible and its closure 
does not disconnect the surface.  Now the condition of  Step $1^\circ$
is that we must make a cut which leaves the set of cuts made
an independent set.
By Theorem~\ref{th:map=Lagrange} maximal independent sets 
({\em bases}\/) have cardinality $n$ and form a 
Lagrangian matroid on $E\cup E^*$. Notice that, due to our definition  
of a basis of a map, independent subsets of (co)edges are just subsets
of bases. 

We use a simplified version of the greedy algorithm for Lagrangian 
matroids as 
described by the following procedure.

Let $\cal B$ be the set of bases of a Lagrangian matroid on $J$. Let
$\cal I$ be the set of all independent subsets, that is all subsets 
of bases from $\cal B$. 
Let $i_1,\ldots, i_n$ be the elements $1,\ldots,n$ written in some order
$i_1 \prec \ldots \prec i_n$.

\bigskip
\noindent
{\sc Greedy Algorithm}

\begin{tabbing}
$B:= \emptyset$;\\
{\tt for $j=1,2,\ldots, n$ do} \\
\hspace{2em} \={\tt if} $B \sqcup \{\,i_j\,\} \in \cal I$ 
{\tt then} $B:=B \sqcup \{\,i_j\,\}$
{\tt else} $B:=B \sqcup \{\,i_j^*\,\}$;\\
{\tt end}
\end{tabbing}

\medskip

Since every basis of a Lagrangian matroid   
contains,  for every $i=1,\ldots,n$,
one of the elements $i,i^*$,
it is immediately obvious that 
the greedy algorithm 
returns a basis in $\cal B$. 

A stronger version of a greedy algorithm, due to Bouchet~\cite{bou},
can be used for characterisation of Lagrangian matroids; 
we do not use it here. 

The greedy algorithm allows us to prove that  the
peeling-off procedure described
produces the set of cuts which is a basis of the map.

Let us pause the procedure at some stage when we have chosen the triangles 
$S_0,S_1,\ldots, S_k$ and let the set of cuts already made be $B$.
Notice  that $B$ is an independent set. 
 If $S_k$ has a cut (co)edge we have to select for $S_{k+1}$
the triangle  lying across the diagonal or non-cut
(co)edge from $S_k$; we can do this unless $S_{k+1}$
is one of the triangles $S_0,\ldots, S_k$;
but in this case one can immediately see that $S_{k+1} = S_0$
and we cut a closed ribbon from the surface.
Since cuts along $B$ do not disconnect the surface,
the ribbon covers the whole surface.

Hence we can assume that neither  edge $e_i$
nor coedge $e_{i^*}$ bounding the
triangle  $S_k$ has not been cut yet. 
We may define an  ordering $\prec$ on $\{\,1,\ldots,n\,\}$
in which non-starred  elements in $B\cup B^*$ preceed all other
elements in $\{\,1,\ldots,n\,\}$.
Obviously, the greedy algorithm for $\prec$ will at some stage produce $B$.
 Furthermore, we may assign $i$ to immediately succeed
 $B$ in the ordering $\prec$, then the greedy
 algorithm says that we can cut along one of
 the edges $e_i$ or $e_{i^*}$ retaining the surface connected,
 and thus continue the process as desired.

Therefore we have proven the following result.

\begin{thm}
Every  simplex in the triangulation of the map $M=(S,G)$
appears in the sequence $S_0,S_1,\ldots, S_n$ exactly once.
The surface $\bar{S}$ obtained by all cuts is homeomorphic to a ring 
{\rm (}punctured $2$-cell\/{\rm )} or\/ $2$-cell.
\end{thm}


\subsection{Representability of Lagrangian matroids} 
\label{subsec:representability}
Let $V$ be the vector space over a field $K$ of characteristic
$\ne 2$ whose basis is
$\{\,e_{i}, e_{i}^* \mid i\in 
I\,\}$. Define a symmetric bilinear form on $V$ by
$\langle  e_{i}, e_{i}^*\rangle =1$,
$ \langle  e_{i}, e_{j}\rangle =0$ for all $i\in I$ and
$j\in J$ with $i^*\neq j$, so that 
the basis $\{\, e_{i}, e_{i}^* \mid i\in I\,\}$
is a {\em hyperbolic basis} in $V$.
$L$ is a {\em Lagrangian} subspace of $V$ if
it is a totally isotropic subspace   
(that is, $\langle k,l\rangle =0$ for every 
$k,l \in L$) of maximal dimension.
If $f_1,\ldots, f_n$ is a basis in a $n$-dimensional
subspace  $L$ in $V$, then we can associate with it
a $n\times 2n$ matrix $C$ written as $(A,B)$ for two $n\times n$ matrices
$A$ and $B$ such that
$$
f_i = \sum_{j=1}^n a_{ij} e_{j} + \sum_{j=1}^n b_{ij} e_j^*.
$$
It is easy to see that $L$ is a Lagrangian subspace in $V$ if and
only if $AB^t$ is a skew symmetric matrix.

Let us index the columns of $A$ with $I$ and those of $B$
with $I^*$, so that the columns of $C$ are indexed by $J$. 
We say that an admissible subset $F\in J$ is {\em independent\/} if
the corresponding columns of $C$ are linearly independent. 
Define the set of {\em bases\/}  
$\B\subseteq J_n$ by putting $X \in \B$ if and only if
\begin{itemize}
 \item $X\in J_n$ and
 \item the $n \times n$ minor consisting of the $i$-th column of $C$ 
for
all $i\in X$
  is non-zero.
\end{itemize}
Notice that change of  a basis in $L$ is equivalent to conducting row
operations on the matrix $(A,B)$ and leaves the pattern of dependencies
of the columns unchanged. Therefore the set 
${\cal B}$ depends only on $L$ and not on choice of basis in $L$.

\begin{thm} {\rm (A.~Vince and N.~White \cite{orth})} \label{absym}
If\/ $L$ is a Lagrangian subspace in $V$  then 
\begin{itemize}
\item every independent subset in $J$ belongs to a basis in $\cal B$, and 
\item $\cal B$ is the set
 of bases of an orthogonal  Lagrangian matroid.
\end{itemize}
\end{thm}

For a proof that the axioms used to define
Lagrangian matroids in \cite{orth} are
equivalent to the Symmetric Exchange Axiom,
see Wenzel \cite{wenz} or the book \cite{book}.

We call an orthogonal matroid arising from a 
Lagrangian subspace $L$ written by a 
matrix $(A,B)$ with $AB^t$ skew-symmetric an {\em orthogonally 
representable} orthogonal matroid, and say that  $(A,B)$ is 
an {\em orthogonal representation} of it over the field 
$K$.\footnote{Note in passing that there is also another way 
of representing Lagrangian matroids, by Lagrangian subspaces 
in a symplectic vector space \cite{bgw}. 
In the more general setting of Coxeter matroids 
\cite{borovik-gelfand} the fact that every orthogonal Lagrangian 
matroid  is a symplectic Lagrangian matroid is explained by the 
canonical embedding of the corresponding 
Coxeter groups $D_n < C_n$ and an observation 
that the  thick $D_n$-building of isotropic subspaces in an 
orthogonal space has the natural structure 
of a thin $C_n$-building \cite{brown}.}

\section{Matroids, Representations and Maps}

Let $\M$ be a map on a compact connected orientable  surface $S$,
with the edge and coedge sets $E = \{\, e_i \mid i \in I\,\}$ and
$E^* = \{\, e_i^* \mid i \in I^*\,\}$, $I = \{\,1,\ldots,n\,\}$,
vertex set $V$ and covertex set $V^*$.
We orient edges $e\in E$ in an arbitrary way 
and choose orientation of coedges $e^*$ so
that the intersection index $(e,e^*) = 1$ for all $e\in E$.

We shall use homology of $S$ and $S \backslash (V \cup V^*)$   
with coefficients in $\mathbb{Q}$.

For a cycle $c$ in $H=H_1 (S \backslash (V \cup V^*))$ we have the
well-defined intersection index, denoted here by $(c,e)$, of 
$c$ with an edge (or coedge) $e \in E \cup E^*$.
We shall denote by $\hat e$ the corresponding linear functional
$\hat{e}: c \mapsto (c,e)$ from $H^*$. Thus $\hat{e}$ is a cocycle
in $H^1(S\setminus (V\cup V^*))$. Notice that this 
$H^1(S\setminus (V\cup V^*))$ can be identified with $H_1(S, V\cup V^*)$
by Poincare-Lefschetz duality \cite[VIII.7]{dold}.

\begin{lem}
An admissible set of {\rm (}co{\rm )}edges $F \subset E \cup E^*$ 
is independent if and only if the linear functionals $\hat{f}$, $f\in F$, 
are linearly independent over $\mathbb{Q}$.
\label{lm:connect=independ}
\end{lem}

\begin{proof}
Denote  $X = V\cup V^* \cup \bigcup_{f\in F} f$. 
Then we have the triple $W \subset X \subset S$ of cell complexes.  
The lemma immediately follows from the long exact homological sequence for
this triple:
\begin{eqnarray*}
\cdots & \longrightarrow &
H_2(X, V\cup V^*) \stackrel{\alpha_2}{\longrightarrow} H_2(S,V\cup V^*)
 \stackrel{\beta_2}{\longrightarrow} H_2(S,X) \\
&&\stackrel{\partial_2}{\longrightarrow} H_1(X, V\cup V^*) 
 \stackrel{\alpha_1}{\longrightarrow} H_1(S,V\cup V^*) \longrightarrow \cdots 
\end{eqnarray*}
Since $H_2(X, V\cup V^*)=0$, $\beta_2$ is an injection.
If $S\setminus X$ is connected then $H_2(S,X)$ is 1-dimensional, 
and, from comparing
dimensions, we see that $\beta_2$ is a surjection. 
Hence $\partial_2=0$ and $\alpha_1$ is an injection.
We need to notice only that $H_1(X, V\cup V^*)$ is generated by (co)edges
in $F$.
\end{proof}

\medskip

Introduce the vector space
${\mathbb{Q}}^E \oplus {\mathbb{Q}}^{E^*}$  over  
$\mathbb{Q}$ with the  basis  $E \cup E^*$.

Define a symmetric bilinear form on ${\mathbb{Q}}^E \oplus {\mathbb{Q}}^{E^*}$ by
$\langle e, e^*\rangle =1$, and $ \langle e,f \rangle =0$ 
for all $e,f \in E\cup E^*$ such that $e^* \ne f$, so that 
 $E\cup E^*$
is a  hyperbolic basis in ${\mathbb{Q}}^E \oplus {\mathbb{Q}}^{E^*}$. 

For $c\in H$, the {\em incidence vector} of $c$, denoted 
$\iota(c)\in {\mathbb{Q}}^E \oplus {\mathbb{Q}}^{E^*}$, is defined by
$$ \iota(c) = \sum_{e \in E\cup E^*} (c,e) e.$$

The main technical result of the paper is:
\begin{thm} \label{cyciso}
    The image  $\iota(H)$ of $H =H_1 (S \backslash (V \cup V^*))$
under the map $c \mapsto \iota(c)$ is an 
    isotropic subspace of the orthogonal
space ${\mathbb{Q}}^E \oplus {\mathbb{Q}}^{E^*}$.
\end{thm}

We postpone the proof of Theorem~\ref{cyciso} until the next section,
and meanwhile deduce from it the main result of the paper,
Theorem~\ref{th:map=Lagrange}.

\begin{thm}\label{goodrep}
The isotropic subspace $\iota(H)$ is Lagrangian, and
the Lagrangian orthogonal matroid corresponding to $\iota(H)$ has 
bases which correspond to the bases of the map $\M$.
\end{thm}

\begin{proof} Notice first that the map $\M = (G,S)$ has at least one basis 
of cardinality $n$.
To construct it, take a spanning tree $T$ in the graph $G$
($T$ can be empty),
then $B=T \cup \overline{T}^*$ is rather obviously a basis of
$\cal M$ and has cardinality $n$
\cite{bou:map}.
If now $b_1,\ldots,b_n$ are the (co)edges in $B$
then the linear functionals
$\hat{B}_1,\ldots, \hat{b}_n$ on $H$ are linearly independent,
hence the $n$ functionals $\tilde{b}_i$ on $\iota(H)$  defined by the rule
$$
\tilde{b}_i(\iota(c)) = \hat{b}_i(c)
$$
are also linearly independent.
Therefore $\dim \iota(H) \ge n$, and being an isotropic
subspace in ${\mathbb{Q}}^E \oplus {\mathbb{Q}}^{E^*}$,
it is a Lagrangian subspace.
Notice that if we represent $\iota(H)$ by a
$n\times 2n$-matrix $M$ in the basis
$e_1,\ldots, e_n, e^*_1,\ldots e^*_n$,
the columns of $M$ will represent the functionals
$\tilde{e}_i, \tilde{e}^*_i$.
By Lemma~\ref{lm:connect=independ} the columns of  $M$
are linearly independent if and only if
the correspondent set of (co)edges is independent in $\cal M$.
Hence the bases of the Lagrangian matroid associated with $\iota(H)$ are
exactly the bases of the map $\M$.
\end{proof}

Now Theorem ~\ref{th:map=Lagrange} is an immediate corollary 
of Theorem~\ref{goodrep}.

\section{\relax Proof of Theorem~\ref{cyciso}}.

We shall prove Theorem~\ref{cyciso} by a naive geometric argument: 
we shall gradually simplify the
map $M$ and use  induction on the total number of vertices and 
covertices in the map. David Stone~\cite{stone} offered a 
cohomological proof of the same result.  

\paragraph{{\sc Reduction, step 1.}}
Notice that if  an edge
$e \in E$ is a contractible loop, i.e.\ the endpoints of
$e$ coincide and the closure $cl(e)$ is a contractible cycle on
$S$ then $(c,e) = 0$ for every $c \in H$ and thus the contribution
$$
(c,e)(d,e^*) + (c, e^*)(d,e)
$$ 
of $e$-th and $e^*$-th coordinates of the incidence vectors
$\iota(c)$ and $\iota(d)$ 
to their scalar product $\langle \iota(c),\iota(d)\rangle $ is $0$.
Hence we can assume  without loss of generality that
the graph $G$ (and, analogously, the dual graph $G^*$)
contains no contractible loops.

\paragraph{{\sc Reduction, step 2.}}  
Introduce the scalar product $\langle\, ,\,\rangle $
on $H$ by setting
$$
\langle c,d\rangle \, = \, \langle \iota(c),\iota(d)\rangle .
$$
Our aim is to prove that $\langle\, ,\,\rangle $ vanishes on $H$.

Consider the kernel $K$ of the canonical map
$$
H_1(S \setminus (V\cup V^*)) \longrightarrow H_1(S).
$$
It is well-known that $K$ is isomorphic to $H_0(V\cup V^*)$ and  generated
by small cycles around vertices and covertices in $V\cup V^*$.

\begin{lem}  \label{lm:(K,H)=0}
$\langle K, H\rangle  \, = 0$.
\end{lem}

\begin{proof}
Indeed, let $c$ be a small cycle around a covertex $v^*\in V^*$ 
(for vertices the proof is analogous).
Let  $e^*_{1,}\ldots, e^*_k$  be 
the oriented edges exiting from or entering  $v$; 
notice that we count loops with both ends at $v^*$ twice,
but with  opposite signs.
The corresponding edges $e_{1},\ldots, e_{k}$ form
the boundary of the country $C$ around $v$; again,
edges corresponding to loops $e^*_i$ appear twice,
with opposite orientations.
Now  if $d$  is an arbitrary cycle in $H$ then 
$$
\langle c,d\rangle \, = \, \langle \iota(c),\iota(d)\rangle \, =
(c,e^*_1)(d, e_{1}) +\cdots + (c,e^*_k)(d, e_{k}).
$$
Without loss of generality we can change  the orientation of
edges and co-edges in such a way that
$(c,e^*_i) = 1$ for all coedges $e^*_i$, $i = 1,\ldots, k$
which are not loops and $e_{1}+\cdots+e_{k}$ is the boundary cycle
$\partial C$ of the country $C$  around $v^*$.
Then the right hand side of the previous equation is exactly
the index of intersection of $d$ and $\partial C$, and hence equals $0$.
\end{proof}

Lemma~\ref{lm:(K,H)=0} allows us to transfer the scalar product 
$\langle\, ,\,\rangle $ from
$H_1(S \setminus (V\cup V^*))$ to $H_1(S)$:
if $c,d \in H_1(S)$ and $c'$, $d'$ are any preimages of
$c$ and $d$ in $H_1(S \setminus (V\cup V^*))$,
then we set $\langle c,d\rangle  \, = \langle  c', d'\rangle $.
The theorem which we are proving  amounts to saying that
$\langle c,d\rangle  \, = 0 $ for all $c,d \in H_1(S)$.
We can work with cycles in $H_1(S)$ instead of
$H_1(S \setminus (V\cup V^*))$ which
gives us the necessary degree of flexibility in
geometric construction we shall use in the proof.

\medskip
 
Now we can start the induction on $|V\cap V^*|$.

\paragraph{\sc {Basis of induction.}}
For the basis of induction assume that the map
$\M$ and the dual map $\M^*$ both contain exactly one vertex,
$v$ and $v^*$, correspondingly.
The orientable connected compact surface $S$ is a sphere with $m$ handles.
Then there are $n=2m$ edges and $2m$ coedges, from the Euler 
formula. We can treat 
the edges  
$e_i$ and coedges $e_i^*$, $i\in I=\{\,1,\ldots n\,\}$  of the map 
as elements of  $H_1(S)$. 
Consider again the kernel $K$ of the canonical surjective map
$$
H_1(S\setminus (V\cup V^*)) \longrightarrow H_1(S);
$$
$K$ is generated by two cycles $k_1$ and $k_2$ around $v$ and $v^*$. 
It is easy to see that  $\iota(k_1) = \iota(k_2)=0$, hence 
$K \leq \ker \iota$ 
and $\iota$ can be lifted to a linear map from $H_1(S)$ to 
${\mathbb{Q}}^E \oplus {\mathbb{Q}}^{E^*}$ which
we denote by the same symbol $\iota$.  Thus 
$\langle c,d\rangle \, = \, \langle \iota(c),\iota(d)\rangle $ 
for $c,d\in H_1(S)$ and 
in our computations we can work with elements of $H_1(S)$ 
instead of $H_1(S\setminus (V\cup V^*))$.

The intersection index 
$$
H_1(S) \times H_1(S) \longrightarrow H_0(S) = {\mathbb{Q}}
$$
is a non-degenerate skew symmetric form on $H_1(S)$. 
From the definition of orientation of  co-edges we have 
$$
(e_i,e_j^*) =  \delta_{ij} \hbox{ for } i,j \in I.
$$
Since $\dim H_1(S) = 2m= |I|$, 
it follows from  elementary linear algebra that 
$\{\, e_i \,\}$ and $\{\, e_i^* \,\}$ are dual 
bases in $H_1(S)$.


If $c\in H_1(S)$ then  we may express $c=\sum_{i \in I} c_i e_i $, 
for some coefficients $c_i$. Now 
we have $\iota(c)$ expressed as a sum of
$\iota(e_i)$. So, if we can show that the space generated by the 
$\iota(e_i)$ is totally isotropic, we have completed
the proof. Now, from the definition of an incidence vector, 
$$\iota(e_i) = \sum_{j\in I}\left( (e_i,e_j)  e_{j} + (e_i, e_j^*) 
 e_j^*\right).$$
Now, taking $i,k \in I$, we have, since $(e_i,e_j^*)=\delta_{ij}$,
\begin{eqnarray*}
\langle \iota(e_i),\iota(e_k)\rangle &=& 
\left\langle \left(\sum_{j\in I} (e_i,e_j) 
 e_{j}+ e_i^*\right),
 \left(\sum_{j\in I} (e_k,e_j)  e_{j}+ e_k^*\right)\right\rangle  \\
&=&\left\langle \sum_{j\in I} (e_i,e_j)  e_{j}, e_k^*\right\rangle 
+\left\langle  e_i^*,\sum_{j\in 
I} (e_k,e_j)  e_{j}\right\rangle  \\
&=&(e_i,e_k)+(e_k,e_i)\\
&=& 0.
\end{eqnarray*}
Thus $\iota(H)$ is generated by a set of elements
with zero pairwise scalar products, so the theorem follows.

\paragraph{\sc {The inductive step.}} Assume that we have a map $\M$ with
$|V \cup V^*| > 2$. Then $\M$ has a contractible edge or coedge, say, edge
$f$. Let $\M'$ be a map on the same surface $S$ obtained by contracting  
the edge $e$ ( in the way shown on Figures~1 and 2) $V'$ and 
${V'}^*$ its sets of vertices and covertices. 
We identify $E'$ with $E \setminus \{\,f\,\}$ 
and ${E'}^*$ with $E^* \setminus \{\,f^*\,\}$, 
correspondingly, and ${\mathbb{Q}}^{E'} \oplus {\mathbb{Q}}^{{E'}^*}$ with 
a subspace of  ${\mathbb{Q}}^{E} \oplus {\mathbb{Q}}^{{E}^*}$, so that the
scalar product $\langle \,,\,\rangle$ on  
${\mathbb{Q}}^{E'} \oplus {\mathbb{Q}}^{{E'}^*}$ is the restriction of 
the scalar product on   ${\mathbb{Q}}^{E} \oplus {\mathbb{Q}}^{{E}^*}$.
The map $\M'$ defines the linear map 
$$
\iota': H_1(S \setminus (V'\cup {V'}^*))\longrightarrow 
{\mathbb{Q}}^{E'} \oplus {\mathbb{Q}}^{{E'}^*},
$$
which allows to lift the scalar product from 
${\mathbb{Q}}^{E'} \oplus {\mathbb{Q}}^{{E'}^*}$ to 
scalar products
$\langle\, ,\, \rangle'$ on 
$H_1(S)$ and $H' = H_1(S\setminus (V'\cup {V'}^*))$.

\begin{lem} \label{induct} 
Contracting a contractable edge $f$ 
does not affect the scalar product\/ 
$\langle\,,\,\rangle$ in the following sense: 
given two cycles $c,d \in H$ we can contract $f$ in such way that, 
for the new map $\M'$ and  appropriate
$c',d'\in H'$, 
$$
\langle c,d,\rangle = \langle c',d' \rangle'.
$$
\end{lem}

\begin{proof} Figures~1 and 2 depict the process of 
contraction.  

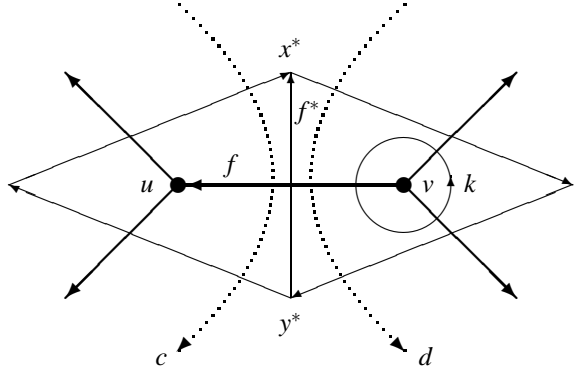
\begin{figure}[h!]

\unitlength 1mm
\begin{picture}(140,60)(-10,0)
\thicklines
\put(60,30){\vector(-1,0){29}}
\put(60,30){\vector(1,1){15}}
\put(60,30){\vector(1,-1){15}}
\put(30,30){\vector(-1,1){15}}
\put(30,30){\vector(-1,-1){15}}
\put(60,8){\vector(1,-1){0.2}}
\bezier{50}(60.00,54.00)(35.00,31.00)(60.00,8.00)  
\put(30,8){\vector(-1,-1){0.2}}
\bezier{50}(30.00,54.00)(55.00,31.00)(30.00,8.00)  

\thinlines
\put(45,45){\line(-5,-2){37.5}}
\put(45,45){\line(5,-2){37.5}}
\put(45,15){\line(5,2){37.5}}
\put(45,15){\line(-5,2){37.5}}
\put(45,15){\vector(0,1){30.00}}
\put(82.5,30){\vector(2,-1){0}}
\put(45,15){\vector(-2,-1){0}}
\put(7.5,30){\vector(-2,1){0}}
\put(45,45){\vector(2,1){0}}
\put(30,30){\circle*{2.00}}
\put(60,30){\circle*{2.00}}
\put(60,30){\circle{12.00}}
\put(66.3,31){\vector(0,1){.2}}
\put(68,29){\makebox(0,0)[bl]{$k$}}
\put(25,29){\makebox(0,0)[bl]{$u$}}
\put(45.5,38){\makebox(0,0)[bl]{$f^*$}}
\put(36,31){\makebox(0,0)[bl]{$f$}}
\put(64,29){\makebox(0,0)[br]{$v$}}
\put(45,47){\makebox(0,0)[b]{$x^*$}}
\put(45,13){\makebox(0,0)[t]{$y^*$}}
\put(27,6){\makebox(0,0)[bl]{$c$}}
\put(62,6){\makebox(0,0)[bl]{$d$}}
\end{picture}
\caption{Fragment of the map $\M$ before contraction.}
\end{figure}

When contracting the edge $f$ we wish to keep 
all changes restricted to the closure of the union of two countries 
of the dual map $\M^*$ containing the endpoints 
$u$ and $v$ of the contracted edge $f$. 
Figure~1 is pre-contraction and Figure~2 is 
post-contraction. The horizontal edge $f= (uv)$ has been 
contracted to a point $u'=v'$ combining the vertices at its ends, 
and the corresponding coedge $f^*=(x^*y^*)$ has disappeared. 
The cycles $c'$ and $d'$ are canonical images of $c$ and $d$
in $H_1(S\setminus (V'\cup {V'}^*))$.

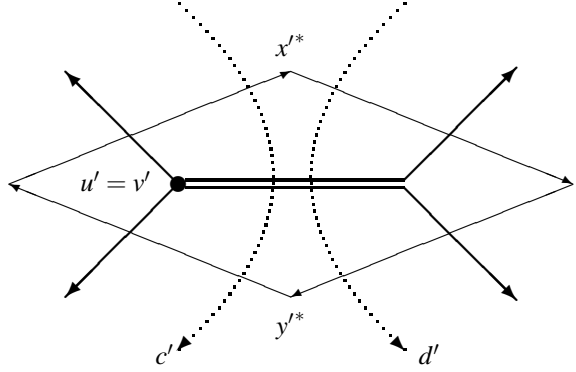
\begin{figure}[h!]
\unitlength 1mm
\begin{picture}(140,60)(-20,0)
\thicklines
\put(60,30.5){\line(-1,0){29}}
\put(60,29.5){\line(-1,0){29}}
\put(60,30.5){\vector(1,1){15}}
\put(60,29.5){\vector(1,-1){15}}
\put(30,30){\vector(-1,1){15}}
\put(30,30){\vector(-1,-1){15}}
\put(60,8){\vector(1,-1){0.2}}
\bezier{50}(60.00,54.00)(35.00,31.00)(60.00,8.00)  
\put(30,8){\vector(-1,-1){0.2}}
\bezier{50}(30.00,54.00)(55.00,31.00)(30.00,8.00)  

\thinlines
\put(45,45){\line(-5,-2){37.5}}
\put(45,45){\line(5,-2){37.5}}
\put(45,15){\line(5,2){37.5}}
\put(45,15){\line(-5,2){37.5}}
\put(82.5,30){\vector(2,-1){0}}
\put(45,15){\vector(-2,-1){0}}
\put(7.5,30){\vector(-2,1){0}}
\put(45,45){\vector(2,1){0}}

%
\put(30,30){\circle*{2.00}}
%
\put(17,29){\makebox(0,0)[bl]{$u'=v'$}}
\put(45,47){\makebox(0,0)[b]{${x'}^*$}}
\put(45,13){\makebox(0,0)[t]{${y'}^*$}}
\put(27,6){\makebox(0,0)[bl]{$c'$}}
\put(62,6){\makebox(0,0)[bl]{$d'$}}
\end{picture}
\caption{Map $\M'$ after contraction.}
\end{figure}

As a result, $c'$, in comparision with $c$, has lost  
its intersections with the coedge $f^*$ and  edge $f$ but has attained  
intersections with every edge  $e\neq f$ which had exited from  $v$ 
prior to contraction; notice, however, that 
the sign of the index of intersection has changed:
$$
(c,e) = -(c',e).
$$
Hence
\begin{eqnarray*}
\iota'(c') & = & \iota(c) -(c,f^*)f^* - (c,f)f 
- (c,f)\cdot \sum_{e \hbox{ {\scriptsize exits from} } v,\; e\ne f} e\\
& = &  \iota(c) -(c,f^*)f^*
- (c,f)\cdot \sum_{e \hbox{ {\scriptsize exits from} } v} e. 
\end{eqnarray*}
If $k$ is a small circle around $v$, as shown on Figure~1, then
$$
\sum_{e \hbox{ {\scriptsize exits from} } v} e\quad = \iota(k)
$$
and 
$$
\iota'(c') =  \iota(c) -(c,f^*)f^* - (c,f)\iota(k). 
$$
Analogously
$$
\iota'(d')  =  \iota(d) -(d,f^*)f^* - (d,f)\iota(k). 
$$
Notice that 
$$
\langle \iota(c),\iota(k)\rangle = \langle \iota(d),\iota(k)\rangle
=\langle\iota(k),\iota(k)\rangle
$$
by Lemma~\ref{lm:(K,H)=0} and $\langle f^*, f^*\rangle=0$ by definition of  
the scalar product $\langle\, ,\,\rangle$ .
Now we can compute:
\begin{eqnarray*}
\langle c',d' \rangle' &= & \langle 
\iota(c) -(c,f^*)f^* - (c,f)\iota(k), 
	\iota(d) -(d,f^*)f^* - (d,f)\iota(k)
\rangle \\
	& = &  \langle \iota(c),\iota(d)\rangle 
	-(d,f^*)\langle \iota(c),f^*\rangle 
		-(c,f^*)\langle f^*, \iota(d)\rangle  \\
	&&\qquad  +(c,f^*)(d,f)\langle f^*,\iota(k) \rangle
				+(c,f)(d,f^*)\langle\iota(k), f^*\rangle.
\end{eqnarray*}
But
$$
\langle \iota(c),f^*\rangle = (c,f), \;\;
\langle f^*, \iota(d)\rangle = (d,f), \;\; 
\langle f^*,\iota(k) \rangle = 1,       
$$
hence we can simplify:
\begin{eqnarray*}
\langle c',d' \rangle' & = &  \langle \iota(c),\iota(d)\rangle  
			     -(d,f^*)(c,f) -(c,f^*)(d,f)\\  
		       && \qquad      +(c,f^*)(d,f)
				+(c,f)(d,f^*)\\
			& = &   \langle \iota(c),\iota(d)\rangle\\
			  & = &   \langle c,d\rangle.
\end{eqnarray*}
\end{proof}

Since contraction of an edge (coedge) does not change the scalar product
$\langle \, , \, \rangle$ on $H_1(S)$, 
we can complete the proof of Theorem~\ref{cyciso} by induction on 
$|V\cup V^*|$.

\section{Non-orientable surfaces}
Obviously in the case when $S$ is a not necessary
orientable compact surface,
we can make all homological computations modulo $2$. 
Symmetric scalar product in that case is also skew symmetric,
and, using obvious 
modifications in terminology and notation, one can prove the following 
two results.

\begin{thm} \label{Thm:cyciso-char2}
    The image  $\iota(H)$ of $H =H_1 (S \backslash (V \cup V^*))$
under the map $c \mapsto \iota(c)$ is an 
    isotropic subspace of the symplectic 
space ${\mathbb{F}}_2^E \oplus {\mathbb{F}}_2^{E^*}$.
\end{thm}

\begin{thm} If\/ $\M$ is a map on a compact surface $S$
then the set\/ $\cal B$ of its bases is a  representable
over ${\mathbb{F}}_2$ symplectic Lagrangian matroid.
\label{th:nonorientedmap=symplectic}
\end{thm}

\subsection*{Acknowledgement}

The authors thank David Stone and Peter Symonds for
valuable discussions.

\end{document}